\documentclass[10pt]{amsart}

\usepackage[colorlinks]{hyperref}
\usepackage{color,graphicx,shortvrb}
\usepackage[latin 1]{inputenc}

\newtheorem{theorem}{Theorem}

\newtheorem{corollary}[theorem]{Corollary}

\newtheorem{definition}[theorem]{Definition}
\newtheorem{lemma}[theorem]{Lemma}

\newtheorem{proposition}[theorem]{Proposition}
\newtheorem{remark}[theorem]{Remark}

\newtheorem{example}[theorem]{Example}
\usepackage{hyperref}

\usepackage[active]{srcltx} 

\usepackage{stackrel}
\usepackage{amssymb}

\def\J#1#2#3{ \left\{ #1,#2,#3 \right\} }
\def\RR{{\mathbb{R}}}
\def\NN{{\mathbb{N}}}
\def\11{\textbf{$1$}}
\def\CC{{\mathbb{C}}}

\usepackage{enumerate}
\usepackage{amsmath}
\usepackage{amssymb}

\begin{document}

\title[Generalised triple homomorphisms]{Generalised triple homomorphisms and Derivations}

\author[Garc{\' e}s]{Jorge J. Garc{\' e}s}
\email{jgarces@ugr.es}
\address{Departamento de An{\'a}lisis Matem{\'a}tico, Facultad de
Ciencias, Universidad de Granada, 18071 Granada, Spain.}

\author[Peralta]{Antonio M. Peralta}
\email{aperalta@ugr.es}
\address{Departamento de An{\'a}lisis Matem{\'a}tico, Facultad de
Ciencias, Universidad de Granada, 18071 Granada, Spain.}

\thanks{Authors partially supported by the Spanish Ministry of Economy and Competitiveness,
D.G.I. project no. MTM2011-23843, and Junta de Andaluc\'{\i}a grants FQM0199 and
FQM3737.}
\thanks{Published at Canadian Journal of Mathematics \url{https://doi.org/https://doi.org/10.4153/CJM-2012-043-7}. This manuscript version is made available under the CC-BY-NC-ND 4.0 license \url{https://creativecommons.org/licenses/by-nc-nd/4.0/} }



\date{}

\begin{abstract} We introduce generalised triple homomorphism between Jordan Banach
triple systems as a concept which extends the notion of generalised homomorphism between
Banach algebras given by Jarosz and Johnson in 1985 and 1987, respectively.
We prove that every generalised triple homomorphism between JB$^*$-triples
is automatically continuous. When particularised to C$^*$-algebras, we rediscover
one of the main theorems established by Johnson. We shall also consider generalised
triple derivations from a Jordan Banach triple $E$ into a Jordan Banach triple $E$-module,
proving that every generalised triple derivation from a JB$^*$-triple $E$ into $E^*$
is automatically continuous.
\end{abstract}
\maketitle
\thispagestyle{empty}

\section{Introduction}

During the last seventy years, a multitude of studies have been published proving that a
homomorphism $T$ between Banach algebras (i.e. a linear map with $T(ab) = T(a)T(b)$
for all $a, b$) must be, under some additional conditions, continuous
(cf. \cite{Dales78}, \cite{Dales00} and \cite{Sin}). For example, it follows
from the original Gelfand's theory that every homomorphism from a Banach algebra to
a commutative, semisimple Banach algebra is automatically continuous. It is well known that
every $^*$-homomorphism between C$^*$-algebras is continuous.
It is due to Johnson that if a unital C$^*$-algebra has no closed cofinite ideals
(e.g. $L(H)$, where $H$ is an infinite dimensional Hilbert space),
then each homomorphism from it into a Banach algebra is continuous (cf. \cite{John}).\smallskip

B.E. Johnson and K. Jarosz considered \emph{generalised homomorphisms} (also called
\emph{$\varepsilon$-multiplicative linear maps} or \emph{$\varepsilon$-isomorphisms}) between Banach algebras
in \cite{Jar85}, \cite{John88} and \cite{John87}. Let $A$ and $B$ be Banach algebras.
A linear mapping $T: A \to B$ is a generalised homomorphism if there exists
$\varepsilon>0$ satisfying $\|T(a b) -T(a) T(b) \| \leq \varepsilon \ \|a\| \ \|b\|$, for every $a,b\in A$.
The first result in this line is due to K. Jarosz, who proved that every generalised homomorphism from a
Banach algebra into a unital abelian C$^*$-algebra is necessarily continuous (cf. \cite[Proposition 5.5]{Jar85}).
B.E. Johnson established in \cite[Theorem 4]{John87} that a generalised homomorphism $T$ between C$^*$-algebras is
continuous if and only if the mapping $a\mapsto T(a^*)^* -T(a)$ is continuous. A \emph{generalised $^*$-homomorphism}
between Banach $^*$-algebras $A$ and $B$ is a generalised homomorphism $T:A \to B$ for
which the mapping $a\mapsto T(a^*)^* -T(a)$ is continuous.\smallskip

Every Banach algebra $A$ can be regarded as an element in the class of Jordan Banach triples
with respect to the product \begin{equation}\label{trivial product Banach algebras}
\J abc := \frac12 (abc +cba).
\end{equation} JB$^*$-triples constitute a subclass of
the Jordan Banach triples which contains all C$^*$-algebras and plays a similar role of that played
by the latter inside the class of Banach algebras (see definitions in section 2).
However, according to our knowledge, the automatic continuity of triple homomorphisms between
Jordan Banach triples (i.e. linear mappings $T$ satisfying $T(\J abc ) = \J {T(a)}{T(b)}{T(c)}$ for every $a,b,c$)
has not been deeply studied. The forerunners in this line reduce to a work of T.J. Barton, T. Dang, and G. Horn,
where these authors prove the automatic continuity of every triple homomorphism between JB$^*$-triples
(see \cite[Lemma 1]{BarDanHor}).\smallskip

In section 3 we define a \emph{generalised triple homomorphism} between Jordan Banach triples $E$ and $F$ as
a linear mapping $T : E \to F$ satisfying $$\|T(\J abc ) - \J {T(a)}{T(b)}{T(c)}\| \leq \varepsilon \ \|a\| \ \|b\| \ \|c\|,$$
for all $a,b,c$ in $E$. We prove that every generalised homomorphism between Banach algebras $A$ and $B$ is
a generalised triple homomorphism when $A$ and $B$ are equipped with the product defined in $(\ref{trivial product Banach algebras})$.
We further show that every generalised $^*$-homomorphism between Banach $^*$-algebras $A$ and $B$ is a generalised triple homomorphism
when $A$ and $B$ are equipped with the product $\J abc := \frac12 (ab^*c +cb^*a)$ (see Proposition \ref{p genhom is gentriple}).
In this section we also study the basic properties of the separating space of a generalised triple homomorphism $T$ between
Jordan Banach triples $E$ and $F$, proving that the separating space $\sigma_{F} (T)$ is a closed triple ideal of the closed
subtriple of $F$ generated by $T(E)$ (compare Proposition \ref{l separating spaces are ideals}).\smallskip

In section 4 we establish some theorems of automatic continuity of generalised triple homomorphism between Jordan Banach triples.
One of the main results proves that every generalised triple homomorphism between JB$^*$-triples is automatically continuous (see Theorem \ref{teo.cont.jbstartriples}). Since every generalised $^*$-homomorphism between C$^*$-algebras is a generalised triple homomorphism,
the aforementioned result of Johnson (see \cite[Theorem 4]{John87}) follows as a direct consequence.
Theorem \ref{t gen.caract.cont} provides necessary and sufficient conditions, in terms of the quadratic annihilator of
the separating space, to characterise when a generalised triple homomorphism from a JB$^*$-triple to a Jordan Banach triple is continuous.
We also prove that every generalised triple homomorphism from a Hilbert space, regarded as a type I Cartan factor, or from a spin factor into
an anisotropic Jordan Banach triple is automatically continuous (cf. Lemmas \ref{l type I Hilbert} and \ref{l spin factor}).\smallskip

In the last section we consider \emph{generalised triple derivations} from a Jordan Banach triple $E$
to a Jordan Banach triple $E$-module $X$.
A conjugate linear mapping $\delta:E\to X$ is said to be a \emph{generalised derivation}
when there exists $\varepsilon>0$ satisfying:
$$\|\delta \{a,b,c\}-\{\delta(a),b,c\}-\{a,\delta(b),c\}-\{a,b,\delta(c)\}\| \leq \varepsilon \ \|a\| \ \|b\| \ \|c\|,$$ for every $a,b, c$ in $E$.
In a recent paper, B. Russo and the second author of this note prove that every triple derivation
from a real or complex JB$^*$-triple $E$ to $E^*$ (i.e., a conjugate linear map $\delta: E \to E^*$ satisfying
$\delta \{a,b,c\}=\{\delta(a),b,c\}+\{a,\delta(b),c\}+\{a,b,\delta(c)\}$)
is automatically continuous (compare \cite[Corollary 15]{PeRu}).
We complement this result proving that every generalised triple derivation from a real or complex JB$^*$-triple
$E$ to $E^*$ is automatically continuous (see Theorem \ref{t genderJB*-triple}). When particularised to C$^*$-algebras,
we show that every generalised triple derivation from a C$^*$-algebra $A$ to a Jordan Banach triple $A$-module
is automatically continuous (compare Theorem \ref{t CtsarToBimodule}). Our results are not mere generalisations
of those forerunners due to Johnson \cite{John87} and Peralta and Russo \cite{PeRu}, the proofs are completely
independent and the theorems presented here are novelties of independent interest even in the category of C$^*$-algebras.

\section{Preliminaries}

We recall that a complex (resp.,
real) \emph{(normed) Jordan triple} is a complex (resp., real)
(normed) space $E$ equipped with a continuous triple product $$ E
\times E \times E \rightarrow E$$
$$(xyz) \mapsto \J xyz $$
which is bilinear and symmetric in the outer variables and
conjugate linear (resp., linear) in the middle one satisfying the
so-called \emph{``Jordan Identity''}:
$$L(a,b) L(x,y) -  L(x,y) L(a,b) = L(L(a,b)x,y) - L(x,L(b,a)y),$$
for all $a,b,x,y$ in $E$, where $L(x,y) z := \J xyz$. If $E$ is
complete with respect to the norm (i.e. if $E$ is a Banach space),
then it is called a complex (resp., real) \emph{Jordan-Banach
triple}. Every normed Jordan triple can be completed in the usual
way to become a Jordan-Banach triple. Unless otherwise is
specified, the term ``normed Jordan triple'' (resp.,
``Jordan-Banach triple'') will always mean a real or complex
normed Jordan triple (resp., ``Jordan-Banach triple'').
\smallskip

For each element $a$ in a Jordan triple $E$, $Q(a)$ will denote
the mapping defined by $Q(a) (x) := \J axa$.\smallskip

Given an element $a$ in a  Jordan triple $E$, we denote $a^{[1]} =
a$, $a^{[3]} = \J aaa$ and $a^{[2 n +1]} := \J a{a^{[2n-1]}}a$
$(\forall n\in \mathbb{N})$. The Jordan identity implies that
$a^{[5]} = \J aa{a^{[3]}}$ , and by induction, $a^{[2n+1]} =
L(a,a)^n (a)$ for all $n\in\NN$. The element $a$ is called
\emph{nilpotent} if $a^{[2n+1]}=0$ for some $n$. Jordan triples
are power associative, that is,
$\J{a^{[k]}}{a^{[l]}}{a^{[m]}}=a^{[k+l+m]}$.\smallskip

A Jordan triple $E$ for which the vanishing of
$\J aaa$ implies that $a$ itself vanishes is said to be
\emph{anisotropic}. It is easy to check that $E$ is anisotropic
if and only if zero is the unique nilpotent element in $E$.\smallskip

A real (resp., complex) \emph{Jordan algebra} is a
(non-necessarily associative) algebra over the real (resp.,
complex) field whose product is abelian and satisfies $(a \circ
b)\circ a^2 = a\circ (b \circ a^2)$. A normed Jordan algebra is a
Jordan algebra $A$ equipped with a norm, $\|.\|$, satisfying $\|
a\circ b\| \leq \|a\| \ \|b\|$, $a,b\in A$. A \emph{Jordan Banach
algebra} is a normed Jordan algebra whose norm is
complete.\smallskip

Every real or complex associative Banach algebra (resp., Jordan
Banach algebra) is a real Jordan-Banach triple with respect to the
product $\J abc = \frac12 (a bc +cba)$ (resp., $\J abc = (a\circ
b) \circ c + (c\circ b) \circ a - (a\circ c) \circ b$).\smallskip

A JB$^*$-algebra is a complex Jordan Banach algebra $A$ equipped
with an algebra involution $^*$ satisfying that $\|\J a{a^*}a \| =
\| 2 (a\circ a^*) \circ a - a^2 \circ a^*\|= \|a\|^3$, $a\in
A$.\smallskip

A \emph{(complex) JB$^*$-triple} is a complex Jordan Banach triple
${E}$ satisfying the following axioms: \begin{enumerate}[($JB^*
1$)] \item For each $a$ in ${E}$ the map $L(a,a)$ is an hermitian
operator on $E$ with non negative spectrum. \item  $\left\|
\{a,a,a\}\right\| =\left\| a\right\| ^3$ for all $a$ in ${A}.$
\end{enumerate}\smallskip

We recall that a subspace $I$ of a normed Jordan triple, $E$, is a
\emph{triple ideal} (resp., a \emph{subtriple}) if
$\{E,E,I\}+\{E,I,E\} \subseteq I$ (resp., $\{I,I,I\} \subseteq I$). The
quotient of a normed Jordan triple by a closed triple ideal is a
normed Jordan triple. It is also known that the quotient of a
JB$^*$-triple (resp., a real JB$^*$-triple) by a closed triple ideal is
a JB$^*$-triple (resp., a real JB$^*$-triple) (compare \cite{Ka}).\smallskip

We recall that a \emph{real JB$^*$-triple} is a norm-closed real
subtriple of a complex JB$^*$-triple (see \cite{IsKaRo95}).\smallskip

A real JB$^*$-algebra is a closed $^*$-invariant real subalgebra
of a (complex) JB$^*$-algebra. Real C$^*$-algebras (i.e., closed
$^*$-invariant real subalgebras of C$^*$-algebras), equipped with
the Jordan product $a\circ b = \frac12 (a b +b a)$, are examples
of real JB$^*$-algebras.

\section{The separating space of a generalised triple homomorphism}

Let $T: E\to F$ be a (non necessarily continuous) linear mapping
between normed Jordan triples. We define $\check{T}: E \times E \times E \to F$ by the rule
$$(a,b,c)\mapsto \check{T}(a,b,c)=T(\{a,b,c\})-\{T(a),T(b),T(c)\}.$$
The mapping $\check{T}$ is symmetric and
linear in the outer variables and conjugate linear in the middle one (trilinear when $E$ is a real
Jordan triple). The mapping $T$ is said to be a \emph{generalised triple homomorphism}
if $\check{T}$ is (jointly) continuous, equivalently, if there exists $C>0$ such that
$$ \|\check{T}(a,b,c)\|=\| T(\{a,b,c\})-\{T(a),T(b),T(c)\}\| \leq C\ \|a\| \|b\| \|c\|. $$

Let $A,B$ be Banach algebras. We have already mentioned that a linear mapping $T:A\to B$ is a
generalised homomorphism  when the bilinear mapping  $$(a,b)\to T(ab)-T(a)T(b)$$ is continuous.
Every Banach algebra is a Jordan Banach triple when endowed with the triple product
\begin{equation}
\label{f trivialTP}2\{a,b,c\}=abc+cba.
\end{equation} We shall refer to this product as the \emph{elemental} (Jordan) triple product of $A.$\smallskip

A richer structure on the Banach algebra $A$ provides richer ternary products. For example,
when $A$ is a $^*$-algebra we can consider the Jordan triple product given by
\begin{equation}
\label{f jordanTP}2\{a,b,c\}=ab^*c+cb^*a.
\end{equation}

Let $A$ and $B$ be Banach *-algebras. A linear mapping $T:A \to B$ is
said to be a \emph{generalised $^*$-homomorphism} if $T$ is a generalised homomorphism
and the mapping $$a\mapsto S(a)=T(a^*)^*-T(a)$$ is continuous. Generalised $^*$-homomorphisms were already considered by B.E.
Johnson in \cite[Theorem 4]{John87}.\smallskip

Our next result explore the connections between generalised (*-)homomorphisms
\hyphenation{homo-morphisms}
and generalised triple homomorphism between Banach (*-)algebras.

\begin{proposition}\label{p genhom is gentriple} Let $A,B$ be Banach
algebras. Every generalised homomorphism $T:A \to B$ is a generalised
triple homomorphism when $A$ and $B$ are equipped with the elemental triple
product $2\{a,b,c\}=abc+cba$.\smallskip

When $A$ and $B$ are
Banach $^*$-algebras and $T$ is a generalised $^*$-homomorphism
then $T$ is a generalised triple homomorphism with respect to the
triple product $2\{a,b,c\}=ab^*c+cb^*a$.

\end{proposition}

\begin{proof}
We start proving the first statement. Let $T:A \to B$ be a generalised
homomorphism between Banach algebras. We shall show that $T$ is a generalised
triple homomorphism when $A$ and $B$ are equipped with the triple
product $(\ref{f trivialTP})$.\smallskip

Throughout this proof, $\tilde{T}$ will denote the continuous bilinear mapping from $A\times A$ into $B$ defined by
$\tilde{T} (a,b) := T(ab)-T(a) T(b).$\smallskip

First, let us see that the (real) trilinear mapping $(a,b,c)\mapsto T(a)\tilde{T}(b,c) $
is continuous. Applying the uniform boundedness principle we see that a trilinear mapping from
the cartesian product of three Banach spaces to another Banach space is (jointly) continuous if, and only if,
it is continuous whenever we fix two variables. Since $\tilde{T}$ is continuous, the desired statement will
follow as soon as we prove that the linear mapping $x\mapsto T(x)\tilde{T}(b,c)$ is continuous whenever we fix $b$ and $c$ in $A$.
Let $(x_n)$ be a norm-null sequence in $A,$ then $$\lim_n T(x_n)\tilde{T}(b,c)=\lim_n T(x_n)T(bc)-T(x_n)T(b)T(c) $$
 $$ =\lim_n \tilde{T}(x_n,b)T(c)+\tilde{T}(x_nb,c)-\tilde{T}(x_n,bc)=0,$$ which proves the desired continuity.\smallskip

Now, the identity $$T(abc)-T(a)T(b)T(c)=\tilde{T}(a,bc)+T(a)\tilde{T}(b,c)$$ implies that the assignment
$(a,b,c)\mapsto T(abc)-T(a)T(b)T(c)$ defines a (jointly) continuous trilinear mapping. It follows that
$$(a,b,c)\mapsto T(\{a,b,c\})-\{T(a),T(b),T(c)\}$$ $$=\frac{1}{2}\Big(T(abc)+T(cba)-T(a)T(b)T(c)-T(c)T(b)T(a)\Big)$$
$$=\frac{1}{2}\Big(T(abc)-T(a)T(b)T(c)\Big)+\frac{1}{2}\Big(T(cba)-T(c)T(b)T(a)\Big)$$
is a continuous trilinear mapping, which gives the first statement.\medskip

Let us suppose now that $T$ is a generalised $^*$-homomorphism
between Banach $^*$-algebras $A$ and $B.$ By the first part of the proof,
$T$ is a generalised triple homomorphism when $A$ and $B$ are equipped with
the triple product $(\ref{f trivialTP})$. We have actually shown that the mapping \begin{equation}
\label{eq first part prop 1} (a,b,c) \mapsto T(abc)-T(a)T(b)T(c)
\end{equation}
is continuous. We shall see that $T$ is a generalised triple homomorphism when $A$
and $B$ are endowed with the product defined in $(\ref{f jordanTP})$.\smallskip

Let us write $S(x)=T(x^*)^*-T(x).$ Fix two elements $a,c$ in $A.$
We claim that the (real) linear mapping \begin{equation}\label{eq claim 1 prop 1}
x\mapsto T(ax^*c+cx^*a)-T(a)T(x)^*T(c)-T(c)T(x)^*T(a)
\end{equation} is continuous.
Clearly, it is enough to check that the restriction to $A_{sa}$ is continuous.
Let $x$ be a self-adjoint element in $A,$ then $$ T(axc)-T(a)T(x)^*T(c)=T(axc)-T(a)T(x)T(c)-T(a)T(x)^*T(c)+T(a)T(x)T(c)$$
$$=T(axc)-T(a)T(x)T(c) - T(a)(T(x)^*-T(x))T(c)$$ and hence
$$T(axc+cxa)-T(a)T(x)^*T(c)-T(c)T(x)^*T(c)$$
$$=\Big(T(axc+cxa)-T(a)T(x)T(c)-T(c)T(x)T(a)\Big)-\Big(T(a)S(x)T(c)+T(c)S(x)T(a)\Big),$$
which proves the claim.\smallskip

Now, we fix $a,b$ in $A$ and claim that the linear mapping \begin{equation}\label{eq claim 2 prop 1}
x\mapsto T(ab^*x)-T(a)T(b)^*T(x),
\end{equation} is continuous. To this end, let $(x_n)$
be a norm null sequence in $A$ then, by $(\ref{eq first part prop 1})$,
$$\lim_n T(ab^*x_{n})- T(a)T(b)^*T(x_n) $$ $$= \lim_n \Big(T(ab^*x_{n}) - T(a)T(b^*) T(x_n)\Big)
+ \Big(T(a)T(b^*) T(x_n) - T(a)T(b)^*T(x_n)\Big)$$
$$=\lim_n \Big(T(ab^*x_{n}) - T(a)T(b^*) T(x_n)\Big) $$ $$+ \lim_n T(a)\tilde{T}(x_n^*,b)^*+T(a)T(b)^*S(x_n) -T(a)\tilde{T}(b^*,x_n)-T(a)S(b^* x_n)=0.$$

Similarly, for every $b,c$ in $A$ the linear mapping
\begin{equation}\label{eq claim 3 prop 1} x\mapsto
  T(xb^*c)-T(x)T(b)^*T(c)
\end{equation} is continuous.

Combining $(\ref{eq claim 1 prop 1})$, $(\ref{eq claim 2 prop 1})$ and $(\ref{eq claim 3 prop 1})$ with the uniform bounded principle
we deduce that the (real) trilinear mapping $(x,y,x)\mapsto T(xy^*z)-T(x)T(y)^*T(z)$ is jointly continuos, and hence, $T$ is a
generalised triple homomorphism for the product defined in $(2)$.
\end{proof}

The separating space of a linear mapping has played and important role in many problems of automatic continuity
(compare, \cite{Rick50}, \cite{Yood}, \cite{BaCur}, \cite{CLev}, \cite{FerGarPer} and \cite{PeRu}, among others).
Let $T: X\to Y$ be a linear mapping between two normed spaces. We recall that the
\emph{separating space}, $\sigma_Y (T)$, of $T$ in $Y$ is defined as the set of all $z$ in $Y$ for which there exists a sequence
$(x_n) \subseteq X$ with $x_n \rightarrow 0$ and $T(x_n)\rightarrow z$. It is well known that a linear mapping $T$ between two Banach
spaces $X$ and $Y$ is continuous if and only if $\sigma_Y (T) =\{0\}.$\smallskip

When $T: A \to B$ is a generalised homomorphism between Banach algebras $A$ and $B$ and $z\in \sigma_Y (T)$
it is not hard to see that $T(a) z$ and $z T(a)$ lie in $\sigma_Y (T)$, for every $a\in A$. This was actually noticed and applied
by B.E. Johnson to show that the separating space, $\sigma_Y (T),$ of $T$ is a closed two-sided ideal of the closed
subalgebra of $B$ generated by $T(A)$ (compare \cite[Lemma 1]{John87}).\smallskip

We are interested in the properties of the separating space of a generalised triple homomorphism $T$ between Jordan Banach triples $E$ and $F$.
Clearly, the image of a generalised triple homomorphism $T: E\to F$ and the image of $\check{T}$ are both contained in
the subtriple of $F$ generated by $T(E).$ However, $T(E)$ and $\check{T} (E\times E\times E)$
need not be Jordan subtriples of $F$. Moreover, it is not so easy to check that the separating space of $T$ is a closed triple ideal of the closed
subtriple of $F$ generated by the image of $T$. The difficulties in the triple setting grow seriously.
For this reason, we shall require an appropriate description
of the subtriple of $F$ generated by a subset. To this end, we define \emph{odd triple monomials}
as follows:\smallskip

Let $E$ be a Jordan triple. An \emph{odd triple monomial of degree $1$ on $E$} is defined as the identity mapping
on $E.$ An \emph{odd triple monomial of degree 3 on $E$} is a mapping $$V:E\times E\times E\to E$$ such that there exists a permutation $\sigma$
of $\{1,2,3\}$ satisfying $$ V(a_1,a_2,a_3)= \{a_{\sigma(1)},a_{\sigma(2)},a_{\sigma(3)} \}.$$
Given a natural $m$, an \emph{odd triple monomial of degree $2m+1$ on $E$} is a mapping $V:E^{2m+1}\to E$
for which there exist a permutation $\sigma$ of $\{1,\ldots,2m+1\}$ and odd triple monomials $V_1,V_2,V_3$ such that
$2m+1=deg(V)=deg(V_1)+deg(V_2)+deg(V_3)$ and $V(a_1,\ldots,a_{2m+1})$ coincides with $$=\{ V_1(a_{\sigma(1)},\ldots,a_{\sigma(i)}),V_2 (a_{\sigma(i+1)},\ldots,a_{\sigma(j)}), V_3 (a_{\sigma(j+1)},\ldots,a_{\sigma(2m+1)})\},$$ where $i=deg(V_1)$ and  $j=deg(V_1)+deg(V_2).$ The odd triple monomial $V$ will be denoted by a pair $V=(\{V_1,V_2,V_3\},\sigma),$ and we shall simply write $V=\{V_1,V_2,V_3\}$ when $\sigma$ is the identity.\smallskip

The symbol $\mathcal{OP}^{^{2m+1}} (E)$ will stand for the set of all odd triple monomials of degree $2m+1$ on $E$, while
$\mathcal{OP} (E)$ will denote the set of all odd triple monomials on $E.$ It should be noticed here that when $V$ is an element in
$\mathcal{OP}^{^{2m+1}} (E)$ and $F$ is another Jordan triple, the mapping $V$ can be regarded as an element in $\mathcal{OP}^{^{2m+1}} (F)$ by
just replacing, in the definition of $V$, the triple product of $E$ with the corresponding triple product on $F$. In order to simplify notation, we shall frequently write $\{ V_1(a_{\sigma(i)}),V_2(a_{\sigma(j)}),V_3(a_{\sigma(k)})\}$ instead of $$\{ V_1(a_{\sigma(1)},\ldots,a_{\sigma(i)}),V_2 (a_{\sigma(i+1)},\ldots,a_{\sigma(j)}), V_3 (a_{\sigma(j+1)},\ldots,a_{\sigma(2m+1)})\}.$$

\begin{lemma}\label{p limpol}
Let $T:E \to F$ be a generalised triple homomorphism between normed
Jordan triples and $m$ a natural number. Let $V$ be an odd triple monomial of degree $2 m+1$ which
can be regarded as an element in $\mathcal{OP}^{^{2m+1}} (E)$ or in $\mathcal{OP}^{^{2m+1}} (F)$.
Suppose $V$ of the form $V=(\{.,W,P\},\sigma)$ $($resp., $V=(\{W,.,P\},\sigma))$, and let $j=deg(W').$
Then $$\lim_{n\to \infty} V(T(x_n),T(a_1),\ldots,T(a_{2m}))- T(V(x_n,a_1,\ldots,a_{2m}))=0,$$
$$\Big(\hbox{resp., } \lim_{n\to \infty} V(T(a_1),\ldots, T(a_{j}),T(x_n), T(a_{j+1}), \ldots, T(a_{2m}))$$
$$- T(V(a_1,\ldots,a_{j},x_n , a_{j+1}, \ldots ,a_{2m}))=0 \Big),$$
for every norm-null sequence $(x_n)$ and $a_1,\ldots,a_{2m}$ in $E.$
\end{lemma}

\begin{proof}
Let $V =(\{V_1,V_2,V_3\},\sigma)$ be an odd triple monomial of degree $2m+1.$
Defining $W=\{V_1,V_2,V_3\}$ then $V(a_1,\ldots,a_{2m+1})=W(a_{\sigma(1)},\ldots,a_{\sigma(2m+1)}).$
Thus, we may assume without loss of generality that $\sigma$ is the identity.\smallskip

We shall proceed by induction on $m$. Since $T$ is a generalised triple homomorphism,
the statement trivially holds for every odd triple monomial of degree $3.$ Now, let us suppose that the
statement is true for odd triple monomials of degree less or equal than $2m-1.$\smallskip

Let $V$ be a odd triple monomial of degree $2m+1.$ We shall assume $V=\{.,W,P\},$ the case $V=\{W,.,P\}$ follows similarly.
Pick a norm-null sequence $(x_n)$ and $a_1,\ldots,a_{2m}$ in $E.$ As we have already observed, we can suppose that
the permutation associated to $V$ is the identity, that is,  $V(x_n,a_1,\ldots,a_{2m})= \{x_n,W(a_i),P(a_j)\}$
$=\{x_n, W(a_1,\ldots,a_{deg(W)}), P(a_{deg(W)+1},\ldots,a_{2m})\}.$\smallskip

The odd triple monomials $W$ and $P$ can be written in the form
$W=\{W_1,W_2,W_3\}$ and $P=\{P_1,P_2,P_3\}$ for some odd triple monomials
$P_i,W_i,$ $i=1,2,3.$ Clearly  $1\leq deg(W_i),deg(P_i)< 2m-1$
(note that we are also assuming that the permutations associated to $W$ and $P$ coincide with the identity
in $\{1,\ldots,deg(V)\}$ and $\{deg(W)+1,\ldots,2m\},$ respectively).\smallskip

Applying the Jordan identity we have {\small \begin{equation}\label{eq -1 lemma 1}V(T(x_n),T(a_1),\ldots,T(a_{2m}))=\{T(x_n),W(T(a_{i})),P(T(a_{j}))\}
\end{equation} $$=\Big\{T(x_n),\{W_1(T(a_{i_1})),W_2(T(a_{i_2})),W_3(T(a_{i_3}))\},\{P_1(T(a_{j_1})),P_2(T(a_{j_2})),P_3(T(a_{j_3})) \}\Big\}$$ $$=\Big\{\Big\{T(x_n),\{W_1(T(a_{i_1})),W_2(T(a_{i_2})),W_3(T(a_{i_3}))\},P_1(T(a_{j_1})) \Big\},P_2(T(a_{j_2})) ,P_3(T(a_{j_3})) \Big\}$$
$$-\Big\{ P_1(T(a_{j_1})), \Big\{\{W_1(T(a_{i_1})),W_2(T(a_{i_2})),W_3(T(a_{i_3}))\},T(x_n),P_2(T(a_{j_2})) \Big\},P_3(T(a_{j_3})) \} $$
$$+ \Big\{ P_1(T(a_{j_1})),P_2(T(a_{j_2})), \Big\{T(x_n), \{W_1(T(a_{i_1})),W_2(T(a_{i_2})),W_3(T(a_{i_3}))\},P_3(T(a_{j_3})) \Big\}\Big\}.$$}\smallskip

We shall treat the summands in the right hand side independently. We claim that {\small
\begin{equation}\label{eq 0 Lemma 1} \lim_n \Big\{\Big\{T(x_n),\{W_1(T(a_{i_1})),W_2(T(a_{i_2})),W_3(T(a_{i_3}))\},P_1(T(a_{j_1})) \Big\},P_2(T(a_{j_2})) ,P_3(T(a_{j_3})) \Big\}
\end{equation}
$$- T\left(\Big\{ \Big\{x_n, \{W_1(a_{i_1}),W_2(a_{i_2}),W_3(a_{i_3}) \},P_1(a_{j_1})\Big\},P_2(a_{j_2}),P_3(a_{j_3}) \Big\}\right)=0.$$}
Indeed, consider the monomial $Q=\{\{., \{W_1,W_2,W_3 \},P_1 \},P_2,P_3 \}.$ It is clear that $deg(Q)\leq 2m-1,$ and
\begin{equation}\label{eq 1 lemma 1} \Big\{T(x_n),\{W_1(T(a_{i_1})),W_2(T(a_{i_2})),W_3(T(a_{i_3}))\},P_1T((a_{j_1})) \Big\}
\end{equation}
$$=Q\Big(T(x_n),T(a_{i_1}),T(a_{i_2}),T(a_{i_3}),T(a_{j_1})\Big).$$

Taking limits in $n\to \infty$ and applying the induction hypothesis we get
\begin{equation}\label{eq 2 lemma 1}
\lim_n Q(T(x_n),T(a_{i_1}),T(a_{i_2}),T(a_{i_3}),T(a_{j_1}) )- T(Q(x_n,a_{i_1},a_{i_2},a_{i_3},a_{j_1}) )=0.
\end{equation}
Let $z_n:=Q(x_n,a_{i_1},a_{i_2},a_{i_3},a_{j_1}).$ It follows from the continuity of the triple product that
$(z_n)$ is a norm-null sequence in $E.$\smallskip

Consider now the monomial $Q'=\{ .,P_2,P_3\}. $  Since $deg(Q')\leq 2m-1$ we can
apply the induction hypothesis to prove
\begin{equation}\label{eq 3 lemma 1}\lim_n \{T(z_n),P_2(T(a_{j_2})),P_3(T(a_{j_3}))\}- T(\{z_n,P_2(a_{j_2}),P_3(a_{j_3})\})
\end{equation}
$$=\lim_n Q'\Big(T(z_n),T((a_{j_2}),T(a_{j_3}))\Big)- T\Big(Q'(z_n,a_{j_2},a_{j_3} )\Big)=0.$$

Combining $(\ref{eq 1 lemma 1})$, $(\ref{eq 2 lemma 1})$ and $(\ref{eq 3 lemma 1})$ we have {\small $$\lim_n
\Big\{\Big\{T(x_n),\{W_1(T(a_{i_1})),W_2(T(a_{i_2})),W_3(T(a_{i_3}))\},P_1(T(a_{j_1}))\Big\},P_2(T(a_{j_2})),P_3(T(a_{j_3}))
\Big\}$$
$$-  T\Big(\Big\{ \Big\{x_n, \{W_1(a_{i_1}),W_2(a_{i_2}),W_3(a_{i_3}) \},P_1(a_{j_1})\Big\},P_2(a_{j_2}),P_3(a_{j_3}) \Big\}\Big)$$
$$=\lim_n
\Big\{Q(T(x_n),T(a_{i_1}),T(a_{i_2}),T(a_{i_3}),T(a_{j_1})),P_2(T(a_{j_2})),P_3(T(a_{j_3}))\Big\}$$
$$-  T\Big(\Big\{ Q(x_n,a_{i_1},a_{i_2},a_{i_3},a_{j_1}),P_2(a_{j_2}),P_3(a_{j_3}) \Big\}\Big)$$
$$=\lim_n \Big\{T(z_n) , P_2(T(a_{j_1})),P_3(T(a_{j_2}))\Big\} - T\Big( \{z_n,P_2(a_{j_1}),P_3(a_{j_2})\}\Big) =0,$$
}
which proves the claim.\smallskip

We can similarly prove that{\small
\begin{equation}\label{eq 4 lemma 1}
\lim_n \Big\{ P_1(T(a_{j_1})), \Big\{\{W_1(T(a_{i_1})),W_2(T(a_{i_2})),W_3(T(a_{i_3}))\},T(x_n),P_2(T(a_{j_2})) \Big\},P_3(T(a_{j_3})) \Big\}
\end{equation}
$$- T\Big(\Big\{ P_1(a_{j_1}), \Big\{\{W_1(a_{i_1}),W_2(a_{i_2}),W_3(a_{i_3})\},T(x_n),P_2(a_{j_2}) \Big\},P_3(a_{j_3})
 \Big\}\Big)=0$$} and
{\small
\begin{equation}\label{eq 5 lemma 1}\lim_n \Big\{\! P_1(T(a_{j_1})),\!P_2(T(a_{j_2})),\! \Big\{T(x_n),\! \{\!W_1(T(a_{i_1})),\!W_2(T(a_{i_2})),\!W_3(T(a_{i_3}))\!\},\!P_3(T(a_{j_3}))\Big\}\Big\}\end{equation}
 $$- T\Big(\Big\{ P_1(a_{j_1}),P_2(a_{j_2}), \Big\{T(x_n), \{W_1(a_{i_1}),W_2(a_{i_2}),W_3(a_{i_3})\},P_3(a_{j_3})
 \Big\}\Big\}\Big) =0.$$}

Finally, from $(\ref{eq -1 lemma 1})$, $(\ref{eq 0 Lemma 1})$, $(\ref{eq 4 lemma 1})$ and $(\ref{eq 5 lemma 1})$ we obtain
{\small
$$\lim_n V(T(x_n),T(a_1),\ldots,T(a_{2m}))- T(V(x_n,a_1,\ldots,a_{2m}))$$
$$= \hbox{ (from the Jordan identity) }=$$
$$\!=\! \lim_n \! \Big\{\!\Big\{\!T(x_n),\!\{W_1(T(a_{i_1})),\!W_2(T(a_{i_2})),\!W_3(T(a_{i_3}))\},\!P_1(T(a_{j_1}))\! \Big\},\!P_2(T(a_{j_2})),\!P_3(T(a_{j_3})) \!\Big\}$$ $$-T\Big(\Big\{ \Big\{x_n, \{W_1(a_{i_1}),W_2(a_{i_2}),W_3(a_{i_3}) \},P_1(a_{j_1})\Big\},P_2(a_{j_2}),P_3(a_{j_3}) \Big\}\Big)$$
$$-\Big\{ P_1(T(a_{j_1})), \Big\{\{W_1(T(a_{i_1})),W_2(T(a_{i_2})),W_3(T(a_{i_3}))\},T(x_n),P_2(T(a_{j_2})) \Big\},P_3(T(a_{j_3})) \Big\}$$
$$+ T\Big(\Big\{ P_1(a_{j_1}), \Big\{\{W_1(a_{i_1}),W_2(a_{i_2}),W_3(a_{i_3})\},T(x_n),P_2(a_{j_2}) \Big\},P_3(a_{j_3})\Big\}\Big)$$
$$+\Big\{ P_1(T(a_{j_1})),P_2(T(a_{j_2})), \Big\{T(x_n), \{W_1(T(a_{i_1})),W_2(T(a_{i_2})),W_3(T(a_{i_3}))\},P_3(T(a_{j_3}))\Big\}\Big\}$$
$$- \lim_n T\Big(\Big\{ P_1(a_{j_1}),P_2(a_{j_2}), \Big\{T(x_n), \{W_1(a_{i_1}),W_2(a_{i_2}),W_3(a_{i_3})\},P_3(a_{j_3})
 \Big\} \Big\}\Big)=0,$$
}as we desired.
\end{proof}

We recall that two elements $a$ and $b$ in a Jordan Banach triple $E$ are
said to be \emph{orthogonal} (written $a\perp b$) if $L(a,b) =
L(b,a)=0$. A direct application of the Jordan identity yields
that, for each $c$ in $E$, \begin{equation}\label{eq orth inner ideal}
a\perp \J bcb  \hbox{ whenever } a\perp b.\end{equation} When $E$ is anisotropic
$a\perp b$ if, and only if, $L(a,b)=0$. In case $E$ is a real or complex JB$^*$-triple,
the relation of being orthogonal admits several equivalent reformulations (cf. \cite[Lemma 1]{BurFerGarMarPe}).
\smallskip

Given a subset $M$ of a Jordan Banach triple, $E,$ we write
$M_{_E}^\perp$ for the \emph{(orthogonal) annihilator of $M$} defined by $$
M_{_E}^\perp:=\{ y \in E : y \perp x , \forall x \in M \}.$$ When
no confusion arise, we shall write $M^{\perp}$ instead of
$M^{\perp}_{_E}$.\smallskip

Let $E$ be a Jordan Banach triple and $S\subseteq E.$ The
norm-closed Jordan subtriple of $E$ generated by $S$ is the
smallest norm-closed subtriple of $E$ containing $S$ and will be denoted by $E_S.$
Clearly, $E_S$ coincides with the norm-closure of the linear span of the
set $$ \mathcal{OP}_{_E} (S) :=\{V(a_{1},\ldots,a_{2m+1}):m\in \NN, V\in \mathcal{OP}^{^{2m+1}} (E),
a_1,\ldots,a_{2m+1} \in S\}.$$ When $a$ in an element in $E$, we write $E_{a}$ instead of $E_{\{a\}}$.

\begin{proposition}\label{l separating spaces are ideals} Let $T: E \to F$ be a
generalised triple homomorphism between two Jordan Banach triples.
Let $I$ denote $\sigma_{F}(T)$ and $F_0$ the norm-closed subtriple of $F$ generated by $T(E)$. Then we have:\begin{enumerate}
\item[$a)$] $I$ is a (closed) triple ideal of $F_0.$
\item[$b)$] $I_{F_0}^{\perp}$ contains all the elements of the
form $\check{T}(a,b,c).$
\end{enumerate}

Further, if $J$ is a closed triple ideal of $F_0$ containing
$I_{F_0}^{\perp},$ then $\pi \circ T$ is a triple homomorphism,
where $\pi$ is the quotient map $F_0 \to F _0/J \cap F_0$.
\end{proposition}

\begin{proof} $a)$ Since $I$ is a closed linear subspace of $F$,
we only have to prove that $\{F_0,F_0,I\}+\{F_0,I,F_0\}\subseteq I.$ Since
$\mathcal{OP}_{_F} (T(E))$ is dense in $F_0$, it is enough to show that
$$V(I,T(E),\ldots,T(E))+V'(T(E),\ldots,T(E),I,T(E),\ldots,T(E))
\subseteq I,$$ where $V$ and $V'$ are arbitrary odd triple monomials of
the form $\{W,.,P\}$ and $\{.,W',P'\},$ respectively.\smallskip

Let $z$ be an element in $I,$ then there exists a
norm-null sequence $(z_n)$ in $E$ such that $z=\lim_n T(z_n).$
Now let $V=\{W,.,P\}$ and $V'=\{.,W',P'\}$ be odd triple monomials
of degree $2m_1+1$ and $2m_2+1,$ respectively, with $j= \hbox{deg} (W)$.
Let us fix $a_1,\ldots,a_{2m_1},$ $b_1,\ldots,b_{2m_2}$ in $E.$
By Lemma \ref{p limpol}

$$ V'(z,T(a_1),\ldots,T(a_{2m_1}))= \lim_n V'(T(z_n),T(a_1),\ldots,a_{2m_1})$$
$$=\lim_n T(V'(z_n,a_1,\ldots a_{2m_1})),$$ and

$$V(T(b_1),\ldots,T(b_{j}),z,T(b_{j+1}),\ldots,T(b_{2m_2}))$$
$$=\lim_n V(T(b_1),\ldots,T(b_{j}),T(z_n),T(b_{j+1}),\ldots,T(b_{2m_2}))$$
$$=\lim_n T(V(b_1,\ldots,b_{j},z_n,b_{j+1},\ldots,b_{2m_2})) .$$

By the continuity of the Jordan triple product  $x_n=V'(z_n,a_1,\ldots a_{2m_1})$ and $y_n=V(b_1,\ldots,b_{j},z_n,b_{j+1},\ldots,b_{2m_2})$
are norm-null sequences in $E$ and thus $$V'(z,T(a_1),\ldots,T(a_{2m_1}))=\lim_n T(x_n) \in I$$ and
$$V(T(b_1),\ldots,T(b_{j}),z,T(b_{j+1}),\ldots,T(b_{2m_2})) =\lim_n T(y_n)\in I.$$\medskip

$b)$ In order to see that $I_{F_0}^{\perp}\supseteq \check{T}(E,E,E),$ we shall show that
$$L(I,\check{T}(a,b,c))|_{F_0}=L(\check{T}(a,b,c),I)|_{F_0}=0, \forall a,b,c \in E.$$

Let $z=\lim T(z_n)$ in $I,$ where $(z_n)$ is a norm-null sequence in $E$, $V$ and odd triple
monomial of degree $2m+1$ and $a,b,c,a_1,\ldots,a_{2m+1}$ in $E.$ Then
$$L(z,\check{T}(a,b,c)) \Big(V(T(a_1),\ldots,T(a_{2m+1}))\Big)$$
$$=\lim_n \Big\{T(z_n),\check{T}(a,b,c),V(T(a_1),\ldots,T(a_{2m+1}))\Big\} $$
$$=\lim_n \{T(z_n),T(\{a,b,c\}),V(T(a_1),\ldots,T(a_{2m+1}))\}$$ $$-\{T(z_n),\{T(a),T(b),T(c)\},V(T(a_1),\ldots,T(a_{2m+1}))\}
=\hbox{ (by Lemma \ref{p limpol}) }=$$ $$=\lim_n T\Big(\!\{z_n,\!\{a,b,c\},\!V(a_1,\ldots,a_{2m+1})\}\Big) -T\Big(\!\{z_n,\!\{a,b,c\},\!V(a_1,\ldots,a_{2m+1})\}\Big)  =0.$$

We can similarly show that $L(\check{T}(a,b,c),z) \Big(V(a_1,\ldots,a_{2m+1})\Big)=0.$ Therefore, it follows from the
density of $\mathcal{OP}_{_F} (T(E))$ in $F_0$ and the continuity of the triple
product that $L(I,\check{T}(a,b,c))|_{F_0}=L(\check{T}(a,b,c),I)|_{F_0}=0,$ which proves $b)$.\smallskip

Finally, to see the last statement we observe that since $I_{F_0}^\perp$ contains all the elements of the
form $\check{T}(a,b,c)$ then we have
$$0=\pi ( \check{T}(a,b,c) )=\pi (T(\{a,b,c\})-\{T(a),T(b),T(c)\})$$
$$=\pi(T(\{a,b,c\}))-\pi(\{T(a),T(b),T(c)\} ),\forall a,b,c \in E,$$
so $\pi \circ T$ is a triple homomorphism. \smallskip
\end{proof}

Let us suppose that, in the hypothesis of Proposition \ref{l separating spaces are ideals} above,
$F$ is assumed to be a JB$^*$-triple. In this setting two elements $a,b$ in $F$ are orthogonal
if and only if $\{a,a,b\}=0$ (cf. \cite[Lemma 1]{BurFerGarMarPe}). Under these assumptions, let $z$ be an element in $I$
and pick arbitrary $a,b,c$ in $E.$ Since there exists a null sequence $(z_n)$ in $E$ such that $z= \lim_{n} T(z_n)$, by Lemma \ref{p limpol},
we have
$$\{z,\!z,\!\check{T}(a,b,c)\} \!=\!\lim_{n} \{T(z_n),\!T(z_n),\!T(\{a,\!b,\!c\}\}-\{T(z_n),\!T(z_n),\!\{T(a),\!T(b),\!T(c)\}\}=0,$$
which implies $I_{F}^{\perp}\supseteq I_{F_0}^{\perp}\supseteq \check{T}(E,E,E)$.

\section{Automatic continuity}

\subsection{Generalised triple homomorphisms between Jordan Banach triples}\ \vspace*{-0.1cm}

A celebrated result of J. Cuntz states that a linear mapping
$T:A\to X$ from a C$^*$-algebra to a Banach space is continuous if,
and only if, its restriction to any C$^*$-subalgebra of $A$ generated
by a single hermitian element is continuous (cf. \cite{Cuntz}). Some years before A.M.
Sinclair \cite{Sin74} established that a similar automatic continuity result holds for
homomorphism from a C$^*$-algebra to a Banach algebra.
At this point, the reader should be tempted to ask if a similar statement holds for linear mappings whose
domain is a JB$^*$-triple (by replacing C$^*$-subalgebras generated by a single hermitian
element by JB$^*$-subtriples generated by a single element). Unfortunately, we shall see next that
the answer to this question is negative.\smallskip

\begin{example}\label{example to Cuntz thm}{\rm
A complex Hilbert space $H$ becomes a JB$^*$-triple when endowed with the triple product defined by
$\{a,b,c\}=\frac{1}{2}((a|b) c + (c|b) a),$ where $(.|.)$ denotes the inner product of $H$. It can be easily seen that
every norm-one element $e$ in $E$ is a tripotent (i.e. $\J eee = e$). Therefore, the JB$^*$-subtriple of $E$ generated
by a single element $a$ coincides with $\CC a$. This implies that, for each Banach space $X$, the restriction of any
linear mapping $T:H\to X$ to any JB$^*$-subtriple of $H$ generated by a single element is continuous. When $H$ is infinite
dimensional we can easily find a discontinuous linear mapping from $H$ into a Banach space.}
\end{example}

The above example shows that a simple translation to the setting of JB$^*$-triples
of the hypotheses assumed by Cuntz in \cite{Cuntz} is not enough to guarantee that a linear mapping from a JB$^*$-triple
to a Banach space is automatically continuous. Finding an assumption to avoid the previous counterexample, we shall replace the
subtriple generated by a single element by the norm-closed inner ideal generated by a single element. We recall that a subspace $J$ of a
JB$^*$-triple $E$ is said to be an \emph{inner ideal} if $\{J,E,J\}$ is contained in $J$. Let $a$ be a norm-one element in $E$
and let $E(a)$ denote the norm closure of $\J aEa$ in $E$. It is known that $E(a)$ coincides with the norm-closed inner
ideal of $E$ generated by $a$ (cf. \cite[Pages 19 and 20]{BuChuZa1}). Let us notice that in the previous Example \ref{example to Cuntz thm},
$H(a) = H$ for every norm-one element $a\in H$.\smallskip

Let $T:E\to F$ be a generalised triple homomorphism between Jordan Banach triples and suppose that $T$
is continuous when restricted to any norm-closed inner ideal generated by a single (norm-one) element.
Let $z$ be an element in $\sigma_F(T).$ Then there exists a norm-null sequence $(z_n)$ in $E$ such that $z=\lim_n T(z_n).$
Pick a norm-one element $a$ in $E.$ Then $$\{T(a),z,T(a)\}=\lim_n \{T(a),T(z_n),T(a)\} = \lim_n T(\{a,z_n,a\}) - \check{T}(a,z_n,a)$$
$$=\lim_n T_{|E(a)} (\{a,z_n,a\}) - \check{T} (a,z_n,a) = 0,$$ since $\{a,z_n,a\}$ is a norm-null sequence in $E(a)$ and
$\check{T}$ and $T_{|E(a)}$ are continuous by hypothesis. Therefore $\{T(z_n),z,T(z_n)\}=0,$ for every natural $n,$ and hence
$z^{[3]}=\lim_n \{T(z_n),z,T(z_n)\}=0,$ which affirms that all elements in $\sigma_F(T)$ are nilpotents.\smallskip

\begin{definition}\label{def CFP} A Jordan Banach triple $E$ has Cohen's factorisation property (CFP)
if given a norm-null sequence $(a_n)$ in $E$ there exist a norm-null sequence $(b_n)$ and
two elements $x,y$ in $E$ such that $a_n=\{x,b_n,y\}, \forall n \in \NN.$
\end{definition}

Every Jordan Banach algebra with a bounded approximate identity has Cohen's factorisation property (compare
\cite{AkLa}). In particular, JB and JB$^*$-algebras have Cohen factorisation property (see \cite[Proposition 3.5.4]{HaOlsen}).
It follows from \cite[Pages 19 and 20]{BuChuZa1} (see also \cite[Lemma 3.2]{EdFerHosPe}) that for every norm-one element
$a$ in a JB$^*$-triple $E$, $E(a)$ satisfies CFP.\smallskip

Our next result is an extension of Sinclair's result \cite[Corollary 4.3]{Sin74}.

\begin{theorem}\label{t gener.inner.ideals}
Let $T:E\to F $ be a linear mapping between two
Jordan Banach triples and suppose that one of the following statements hold:
\begin{enumerate}
 \item[$a)$] $T$ is a generalised triple homomorphism and $F$ is anisotropic;
 \item[$b)$] $E$ has Cohen's factorisation property.
\end{enumerate}
If the restriction of $T$ to any closed inner ideal generated by a
single element is continuous, then $T$ is continuous.
\end{theorem}

\begin{proof}
The proof under hypothesis $a)$ was already given in the paragraph preceding Definition \ref{def CFP}.
Suppose $E$ satisfies CFP. Let $(y_n)$ be a norm-null sequence in $E$ and let $a\in E$.
Since $T|_{E(a)}$ is continuous, we have $\lim_n T\{a,y_n,a\}= 0$. Since $a$ was arbitrary chosen,
we deduce that
\begin{equation}
\label{eq 1 them Sinclair}\lim_{n} T(\{a,y_n,b\})= 0,
\end{equation} for every $a,b\in E$.\smallskip

Let us pick $z\in \sigma_{F} (T)$ and a norm-null sequence $(z_n)$ in $E$
satisfying $T(z_n) \to z$. By hypothesis, there exist
$a,b$ in $E$ and a norm-null sequence $(y_n)\subseteq E$ such that $z_n= \J a{y_n}b$. In such a case,
by $(\ref{eq 1 them Sinclair})$, $$z= \lim_n T(z_n) =  \lim_n T(\J a{y_n}b) = 0.$$
\end{proof}

\begin{remark} \label{sep.sapce.composition}{\rm Let $T:E\to F$ be a linear mapping between Banach spaces.
A useful property of the separating space
$\sigma_F(T)$ asserts that for every bounded linear $R$ from $F$
to another Banach space $Z,$ the composition $RT$ is continuous if, and only if, $\sigma_F(T)\subseteq ker(R).$
It is also known that $ \sigma(RT)=\overline{R(\sigma(T))}^{\|.\|}$ (see Lemma 1.3 in \cite{Sin}.)}
\end{remark}

Based on the Commutative Gelfand Theory established
by W. Kaup (cf. \cite{Ka}), T.J. Barton, T. Dang, and G. Horn proved
the automatic continuity of triple homomorphisms between JB$^*$-triples (see \cite[Lemma 1]{BarDanHor}).
The natural extension of this automatic continuity property to the setting of generalised triple
homomorphisms is contained in our next result.

\begin{theorem}\label{teo.cont.jbstartriples}
Every generalised triple homomorphism between JB$^*$-triples is continuous.
\end{theorem}

\begin{proof} Let $T : E \to F$ be a generalised triple homomorphism between JB$^*$-triples.
The norm closed subtriple of $F$ generated by $T(E)$ will be again denoted by $F_0$, while the symbol
$I$ will stand for the separating space $\sigma_{F} (T)$.
Since $F_0$ is a norm-closed subtriple of $F,$ then $F_0$ is a
JB$^*$-triple itself. Proposition \ref{l separating spaces are ideals} $a)$ assures that $I$ is a closed ideal of $F_0$, and by \cite[Lemma 4]{PeRu} $I_{F_0}^{\perp}$ is a norm-closed triple ideal of $F_0$.\smallskip

The final statement in Proposition \ref{l separating spaces are ideals} guarantees that
the linear mapping $\pi\circ T: E \to F/I_{F_0}^{\perp}$ is a triple homomorphism. Since the quotient
$F_0/I_{F_0}^{\perp}$ is a JB$^*$-triple, the triple homomorphism $\pi \circ T$ is continuous (cf. \cite[Lemma 1]{BarDanHor}).
By remark \ref{sep.sapce.composition} we have $I=\sigma_F(T)\subseteq ker(\pi)=I_{F_0}^{\perp}$ and the latter implies that
$I=\sigma_F(T)=0.$
\end{proof}

Since every $C^*$-algebra endowed with the triple product given in $(\ref{f jordanTP})$,
is a JB$^*$-triple, Theorem \ref{teo.cont.jbstartriples} together with Proposition \ref{p genhom is gentriple}
allow us to rediscover the following result which is originally due to B.E. Johnson \cite [Theorem 4]{John87}.

\begin{corollary}[Theorem 4, \cite{John87}]
Every generalised $^*$-homomorphism between $C^*$-algebras is
continuous. $\hfill \Box$
\end{corollary}

Our next goal is to explore the automatic continuity of a generalised triple homomorphism from
a JB$^*$-triple to a Jordan Banach triple. To this end we shall require some additional concepts
and tools.\smallskip

Let $E$ be a real or complex JB$^*$-triple. We shall
say that E is \emph{algebraic of bounded degree} if all singly generated
subtriples of $E$ are finite-dimensional. If in fact there exists
$m\in \NN$ such that single-generated subtriples of $E$ have
dimension $\leq m,$ then $E$ is said to be of \emph{bounded degree,}
and the minimum os such an $m$ will be called the degree of $E.$\smallskip

Our next result owes much to the proof given in \cite[Proposition 12]{PeRu}
by B. Russo and the second author of this note.

\begin{theorem}\label{t AutoContToBanachSpace}
Let $T: E \to X$ be a linear mapping from a JB$^*$-triple to a
Banach space. Let $J_T:=\{a\in E: T\circ Q(a), T\circ L(a,a)
\mbox{ are continuous} \}.$ Suppose that $J_T$ 
has the following properties:
\begin{enumerate}[$(a)$]
\item $a+b$ lies in $J_T$ whenever $a,b\in J_T$;
\item $\{E,E,J_T\}+\{E,J_T,E\} \subseteq J_T;$
\item If $I$ is a norm-closed triple ideal containing $J_T$ then
$E/{I}$ is algebraic of bounded degree.
\end{enumerate}
Then $T$ is continuous if and only if $J_T$ is norm-closed.\end{theorem}

\begin{proof} When $T$ is continuous, $J_T$ coincides with $E$ and nothing has to be proved.
Suppose now that $J_{T}$ is norm-closed. It follows from $(a)$ and $(c)$ that $J_T$ is a norm-closed triple ideal
of $E$. We claim that the restriction of $T$ to $J_T$ is continuous. Indeed, the assignment $(a,b,c)\mapsto W(a,b,c)=T(\{a,b,c\})$
defines a (real) trilinear mapping $W:J_T\times J_T\times J_T \to F$.  From the
definition of $J_T$ and $a),$ $W$ is separately continuous
whenever we fix two variables. An application of the uniform boundedness principle implies that
$W$ is jointly continuous. Therefore, there
exists a positive constant $M$ such that $\| T \{a,b,c\} \|
\leq M \|a\| \|b\| \|c\|,$ for every $a,b,c$ in $J_T.$
Since $J_T$ is a JB$^*$-subtriple of $E,$ for each $a$ in $J_T$ there
exists $b$ in $J_T$ such that $b^{[3]}=a.$ In this case
$$\|T(a)\|=\|T(\{b,b,b\})\| \leq M \|b\|^3 = M \ \| \{b,b,b\}\|=M\ \|a\|,$$
which shows that $T|_{J_T}$ is continuous.\smallskip

Finally, let us prove that $J_T=E.$ By hypothesis $(c)$,
$E/J_T$ is algebraic of bounded degree $m.$ Thus, for each
element ${a}+J_T$ in $E/J_T$ there exist mutually orthogonal
minimal tripotents ${e_1}+J_T,\ldots,{e_k}+J_T$ in
$E/J_T$ and $0<\lambda_1\leq\ldots\leq \lambda_k$ with $k\leq m$
such that ${a}+J_T=\sum_{j=1}^k \lambda_k {e_k}+J_T.$
We shall show that $e_1,\ldots,e_k \in J_T,$ and hence, $a\in J_T$,
which proves $E=J_T$.\smallskip

Let ${e}+J_T$ be a minimal tripotent in $E/J_T$. Henceforth, $\pi: E\to E/J_T$
will denote the canonical projection. Take an arbitrary norm-null
sequence $(a_n)$ in $E$. For each natural $n,$ there exists a scalar $\mu_n\in \CC$
such that $\pi (Q(e)(a_n))=\mu_n  ({e}+ J_T).$ The continuity of $\pi$ and the Peirce
projection $P_2({e}+J_T)$ assure that $\mu_n \to 0.$ It follows that $Q(e)(a_n)-\mu_n e$
lies in $J_T$ and tends to zero in norm. Since, by hypothesis, $T|_{J_T}$ is continuous we have
$$ T( Q(e)(a_n))=T( Q(e)(a_n)-\mu_n e) +\mu_n T(e) \to 0.$$ The arbitrarity of $(a_n)$ guarantees
that $T\circ Q(e)$ is continuous, or equivalently, $e$ lies in $J_T.$
\end{proof}

The following auxiliar lemmas will be needed later.

\begin{lemma}\label{l bded degree.gral}
Let $E$ be a real JB$^*$-triple and $J$ a subset of $E$ satisfying
that whenever we have two sequences $(x_n),(y_n)$
in $E$ such that $Q(y_n)Q(x_n)=Q(x_n)$ and $Q(y_n)Q(x_m)=0$ for $n\neq m$,
then $x_n$ lie in $J$ except (perhaps) for finitely many $n.$
Suppose $I$ is a norm-closed triple ideal of $E$ containing $J$ then $E/I$
is algebraic of bounded degree.
\end{lemma}

\begin{proof} Since $I$ contains $J$ then $I$ also has the
property assumed in the hypothesis.\smallskip

Let us write $F=E/{I}.$ As noticed in the proof
of Corollary 8 in \cite{PeRu} for $\overline{a}=a+{I}$ we have
$F_{\overline{a}}=E_a/(E_a\cap I).$\smallskip

The commutative JB$^*$-triple $E_a$ is triple isomorphic to some $C_0(L)$ (cf. \cite[\S 1]{Ka}).
We shall identify $E_a$ with $C_0(L).$ It is known that $F_{{a}+I}\cong C_0(\Gamma)$
where $$\Gamma=\{t \in L: b(t)=0, \forall b\in E_a \cap I\}.$$
We claim that $\Gamma$ is finite. Otherwise, there exists an infinite sequence $(t_n)$ in $\Gamma$
and a sequence of open disjoint sets $\{U_n\}_n.$ By local compactness we can find open sets $V_n,W_n$
with $\overline{V_n}$ and $\overline{W_n}$ compact, such that $t_n \in V_n \subseteq
\overline{V_n} \subseteq W_n \subseteq \overline{W_n}\subseteq U_n.$\smallskip

By Urysohn's lemma, for each natural $n$, we can find $f_n\in C_0(L)$ with $t_n\in supp(f_n)\subseteq W_n$ and
$g_n\in C_0(L)$ such that $g_n\equiv 1$ in $\overline{W_n}$ and vanishing outside $U_n.$
Since $f_n(t_n),g_n(t_n)\neq 0, \forall n \in \NN$ then $f_n,g_n \notin I,\forall n\in \NN.$
In this case the sequences $(f_n),(g_n)$ verify that $Q(g_n)Q(f_n)=Q(f_n)$ and $Q(g_n)Q(f_m)=0$ for $n\neq m,$
and they do not lie in $I$, which is a contradiction.\smallskip

It follows that $\Gamma$ is finite and therefore $F_{{a}+I}$ is finite dimensional.
Since $a+I$ was arbitrary chosen, the statement of the lemma follows from \cite[Theorem 3.8]{BeLoPer}.\smallskip
\end{proof}

\begin{lemma}\label{l gen QTcontinua} Let $T:E\to F$ be a generalised triple homomorphism between
real Jordan Banach triples and let $(x_n),(y_n)$ be sequences of elements in $E$ such that $Q(y_n)Q(x_n)=Q(x_n)$
and $Q(y_n)Q(x_m)=0$ for $n\neq m.$ Then $Q(T(x_n))T$ and $T Q(x_n)$ are continuous for all but a finite number of $n$.
 \end{lemma}

 \begin{proof} Let us suppose that $Q(T(x_n))T$ is discontinuous for
 infinitely many $n$ in $\NN.$ By passing to a subsequence if necessary, we can assume that
 $Q(T(x_n)) T$ is discontinuous for all $n$ in $\NN.$ We observe that, since $T$ is a generalised triple homomorphism
 the identity
 $$\{T(x_n),T(b),T(x_n)\} =T(\{x_n,b,x_n\})-\check{T}(x_n,b,x_n),$$ holds for every $b\in E$ and $n\in \mathbb{N}$.
 It is then clear that $Q(T(x_n)) T$ is continuous
 if, and only if, $TQ(x_n)$ is. So, we may assume that $T Q(x_n) $ is
 discontinuous for all $n$ in $\NN$. Choose $(a_n)$ in $E$ such that $\|a_n\|\leq
 2^{-n}\|x_n\|^{-2}$ and $$\|T Q(x_n)(a_n)\|\geq 2^n ( 1+\|T(y_n)\|^2) +
 \| \check{T}\| \|y_n\|^2.$$ Let $a=\sum_{m\geq 1} \{ x_m,a_m,x_m\}.$ Since
 $\{y_n,a,y_n\}=\{x_n,a_n,x_n\}$ we have $$2^n (1+ \|T(y_n)\|^2) +
 \| \check{T}\| \|y_n\|^2 \leq \|T Q(x_n)(a_n)\|$$ $$= \|T Q(y_n)(a)\|=\| Q(T(y_n))(T(a))+\check{T}(y_n,a,y_n)\| $$
$$ \leq \|T(y_n)\|^2 \|T(a)\| +\|\check{T}\|\|y_n\|^2\|a\| \leq (1+\|T(y_n)\|^2)\|T(a)\|+ \|\check{T}\|\|y_n\|^2.
$$ So we have that $\|T(a)\| \geq 2^n, \forall n \in \NN,$ which is impossible.
\end{proof}

Let $T:E\to F$ be a generalised triple homomorphism between Jordan Banach triples.
Following the notation employed in Proposition \ref{l separating spaces are ideals},
the symbol $F_0$ will denote the norm-closed subtriple of $F$ generated by $T(E)$.\smallskip

According to the notation defined in \cite{PeRu}, for each subset $B$ of a Jordan Banach triple $F$,
we define its \emph{quadratic annihilator}, Ann$_{F} (B)$, as the set $$\{ a\in
F : Q (a) (B) = \J aBa = 0\}.$$ The quadratic annihilator will be used later in a more general setting.\smallskip

If we set $J:=T^{-1}(Ann_F(\sigma_F(T))),$ it not hard to see, from the basic properties of the separating space,
that $J$ coincides with the set $\{a\in E: Q(T(a)) T \mbox{ is continous}\}$
(compare Remark \ref{sep.sapce.composition}), and since $T$ is a generalised
triple homomorphism, the latter equals $\{a\in E: T Q(a) \mbox{ is continous}\}$ (compare the proof of Lemma \ref{l gen QTcontinua}). The following result follows straightforwardly from Lemmas \ref{l gen QTcontinua} and \ref{l bded degree.gral} and the above comments.

\begin{proposition}\label{t E/J hilbert} Let $T:E\to F$ be a generalised triple homomorphism from a real
JB$^*$-triple to a Jordan Banach triple. The following statements hold:
 \begin{enumerate}[$a)$]
\item If $I$ is a norm-closed triple ideal containing $T^{-1}(Ann_F(\sigma_F(T)))$ then
$E/{I}$ is algebraic of bounded degree.

\item Let $K$ be a triple ideal of E. The linear mapping
 $$x\in E \mapsto \{T(a),T(x),T(a)\}$$ is continuous for all $a$ in $K$ if,
 and only if, $K$ is contained in $T^{-1}(Ann_F(\sigma_F(T))).$ 
\end{enumerate}
\end{proposition}

We can establish now the main result of this section.

\begin{theorem}\label{t gen.caract.cont} Let $T:E\to F$ be a generalised triple homomorphism
from a JB$^*$-triple to a Jordan Banach triple and let $J=T^{-1}(Ann_F(\sigma_F(T))).$
The following statements are equivalent  \begin{enumerate}[$a)$] \item $J$ is a norm-closed triple ideal of $E$ and
$$\{Ann_F(\sigma_F(T)),Ann_F(\sigma_F(T)),\sigma_F(T)\}=0.$$ \item $T$ is continuous.
\end{enumerate}
\end{theorem}

\begin{proof} The implication $b)\Rightarrow a)$ is clear. We shall prove $a)\Rightarrow b)$.
We already know, by Proposition \ref{t E/J hilbert} $b)$,
that for each element $a$ in $J$, the linear mapping $$x \in
E \mapsto \{T(a),T(x),T(a)\}$$ is continuous. Let us fix
$a,b$ in $J.$ Since $J$ is a linear subspace of $E$ then $a+b$
also lies in $J,$ that is, the mapping $x
\mapsto\{T(a+b),T(x),T(a+b)\}$ is continuous. The identity
$$2\{T(a),T(x),T(b)\}=\{T(a+b),T(x),T(a+b)\}$$ $$-\{T(a),T(x),T(a)\}-\{T(b),T(x),T(b)\},$$
guarantees that the mapping $x\mapsto \{T(a),T(x),T(b)\}$ is continuous, or
equivalently, $T Q(a,b)$ is continuous (because $T$ is a
generalised triple homomorphism).\smallskip

Since $\{Ann_F(\sigma_F(T)),Ann_F(\sigma_F(T)),\sigma_F(T)\}=0,$ the
linear mapping $$x \in E \mapsto \{T(a),T(b),T(x)\}$$ is continuous for every $a,b\in J$. Applying that
$T$ is a generalised triple homomorphism, we deduce that the
linear mapping $x \in E \mapsto T(\{a,b,x\})$ also is continuous for every $a,b\in J$.
This shows that the trilinear mapping $W:E\times E\times E,$ given by
$(a,b,c)\mapsto W(a,b,c)= T(\{a,b,c\})$ is continuous whenever we
fix two variables in $J$. An application of the unform boundedness principle proves that
$W|_{J\times J\times J}$
is jointly continuous. Following the argument given in the proof of Theorem \ref{t AutoContToBanachSpace}, we show that $T|J: J\to F$ is continuous.\smallskip

Proposition \ref{t E/J hilbert} $a)$, implies that $E/J$ is algebraic of bounded degree.
The proof concludes applying the argument given in the final part of the proof
of Theorem \ref{t AutoContToBanachSpace}.
\end{proof}

The above Theorem \ref{t gen.caract.cont} admits a more detailed statement in the
particular setting of some Cartan factors.\smallskip

We recall that a complex Hilbert space $H$ can be regarded as a \emph{type I Cartan factor}
with its natural norm and the product given by $$2 \J abc := (a|b) c + (c|b) a, \ \ (a,b,c\in H),$$ where
$(.|.)$ denotes the inner product of $H$.

\begin{lemma}\label{l type I Hilbert} Let $H$ be a complex Hilbert space regarded as a type I Cartan factor, $F$ an anisotropic
Jordan Banach triple and $T: H \to F$ a generalised triple homomorphism. Then $T$ is continuous.
\end{lemma}

\begin{proof} Let $F_0$ denote the norm-closed subtriple of $F$ generated by $T(E)$. It is enough to prove
that $T: H \to F_0$ is continuous. Replacing $F$ with $F_0$,
we may assume, by Proposition \ref{l separating spaces are ideals}, that $\sigma_{F} (T)$ is
a norm-closed triple ideal of $F$ and $F$ is generated by $T(E)$.
It follows from our hypothesis that the mapping $$(a,b,c)\mapsto \check{T} (a,b,c)=\frac12 \left((a|b) T(c) + (c|b) T(a)\right) - \J {T(a)}{T(b)}{T(c)}, \ \  (a,b,c\in H),$$ is continuous. Let $z$ be an element in $\sigma_{F} (T)$, there exists a norm-null sequence $(x_{n})\subset H$ such that $T(x_n) \to z$.
If we fix two arbitrary elements $a$, $c$ in $H$, it follows from the continuity of $\check{T}$ and the triple product of $F$ that $$0=\lim_{n} \frac12 \left((a|x_n) T(c) + (c|x_n) T(a)\right) - \J {T(a)}{T(x_n)}{T(c)} = - \J {T(a)}{z}{T(c)}.$$ It follows from the arbitrariness of $a$ and $c$ that $\J {T(E)}{\sigma_{F}(T)}{T(E)}=0.$ Similarly, let $V$ and $W$ be odd triple monomials of degree $2m_1+1$ and $2m_2+1,$ respectively, and let us fix
$a_1,\ldots,a_{2m_1},$ $b_1,\ldots,b_{2m_2}$ in $H.$ By Lemma \ref{p limpol}, $$\J {V(T(a_1),\ldots, T(a_{2 m_1+1}))}{z}{W(T(b_1),\ldots, T(b_{2 m_2+1}))}$$ $$= \lim_{n} \J {V(T(a_1),\ldots, T(a_{2 m_1+1}))}{T(x_n)}{W(T(b_1),\ldots, T(b_{2 m_2+1}))}$$ $$= \lim_n  T\left(\J {V(a_1,\ldots, a_{2 m_1+1})}{x_n}{W(b_1,\ldots, b_{2 m_2+1})}\right) $$
$$= \lim_n   \frac12 ({V(a_1,\ldots, a_{2 m_1+1})}|{x_n})\  T\left({W(b_1,\ldots, b_{2 m_2+1})}\right) $$ $$+ \frac12 ({W(b_1,\ldots, b_{2 m_2+1})}|{x_n}) \ T\left({V(a_1,\ldots, a_{2 m_1+1})}\right) =0.$$ Since we have assumed that $F$ is the Jordan Banach triple generated by $T(E)$, it follows by linearity and from the continuity of the product of $F$ that $$\J {F}{\sigma_{F}(T)}{F}=0.$$ Finally, $F$ being anisotropic implies that $\sigma_{F}(T)=0$ and hence $T$ is continuous.
\end{proof}

A \emph{(complex) spin factor} is a complex Hilbert space $S$ provided with
a conjugation (i.e. a conjugate linear isometry of period 2) $x \mapsto \overline{x}$,
triple product $$\J abc = \frac12 \left( \left( a | b \right) c + \left( c | b \right) a - \left( a |
\bar c \right) \bar b\right),$$ and norm given by $\|a\|^2=\frac12\ \left( a | a
\right)+\frac12 \ \sqrt {\left( a | a \right)^2-|\left( a | \overline a
\right)|^2}$, for every $a,b,c\in S$.

\begin{lemma}\label{l spin factor} Let $S$ be a (complex) spin factor, $F$ an anisotropic
Jordan Banach triple and $T: S \to F$ a generalised triple homomorphism. Then $T$ is continuous.
\end{lemma}

\begin{proof} Let $S$ be a spin factor. Corollary in \cite[page 313]{DaFri} and the proof of
Proposition in page 312 in the just quoted paper assures that $S$ is the norm closed linear span of a
``\emph{spin grid}'' $\{u_i,v_i, u_0 \}_{i\in \Gamma}$, where $(u_i|u_j)=(v_i|v_j)
=(u_i|v_j)=(u_i|v_i)= (u_0|u_i)= (u_0|v_i) =0,$ $\|u_i\|=1$, $\|v_i\|=1$, $\|u_0\|=1 \hbox{ or } 0$, $\overline{u_i}= v_i,$
and $\overline{u_0} = u_0,$ for every $i\neq j$ in $\Gamma$. Let $S_1$ (resp., $S_2$) denote the norm-closed subspace of
$S$ generated by $\{u_i: i\in \Gamma\}$ (resp., $\{v_i: i\in \Gamma\}$). Clearly $S = S_1 \oplus S_2 \oplus \mathbb{C} u_0$.
It is easy to see that $S_1$ and $S_2$ are norm-closed subtriples of $S$ (i.e. $\J {S_i}{S_i}{S_i} \subset {S_i}$) and
$\J abc = \frac12 \left( \left( a | b \right) c + \left( c | b \right) a \right),$ for every $a,b,c$ in $S_i$ ($i=1,2$).
Therefore $S_1$ and $S_2$ are Hilbert spaces equipped with structure of type I Cartan factors. Lemma \ref{l type I Hilbert}
shows that $T|_{S_i} : S_i \to F$ is continuous for every $i=1,2$. Finally, the continuity of the natural projections of $S$
onto $S_1$, $S_2$ and $\mathbb{C} u_0$ assures that $T$ is continuous.
\end{proof}

According to the comments given before Proposition 17 in \cite{PeRu}, the proof of Theorem \ref{t AutoContToBanachSpace}
(and hence the proof of Theorem \ref{t gen.caract.cont}) is only valid for complex JB$^*$-triples, the reason being that,
in the real setting, a minimal tripotent $e$ in a real JB$^*$-triple $E$ need not satisfy that $E_2 (e) = \mathbb{R} e$.
Actually, there exist examples of minimal tripotents $e$ for which $E_2 (e)$ is infinite dimensional.
The extension of Theorem \ref{t gen.caract.cont} to the real setting is not a trivial
consequence of the result proved in the complex case and constitute a result of independent interest
which remains open in this paper. However, there exists a subclass of real JB$^*$-triples for which the statements of Theorems \ref{t AutoContToBanachSpace}
and \ref{t gen.caract.cont} remains true. A real JB$^*$-triple $E$ is called \emph{reduced} whenever $E_{2} (e) = \RR e$
(equivalently, $E^{-1} (e) = 0$) for every minimal tripotent $e\in E$.
Reduced real JB$^*$-triples were considered in \cite{Loos77}, \cite{Ka97}, \cite{FerMarPe} and \cite{PeRu}.
The proof of Theorem  \ref{t gen.caract.cont} is valid for reduced real JB$^*$-triples.\smallskip

\subsection{Generalised triple derivations from a JB$^*$-triple}\medskip

B. Russo and the second author of this note carried out in \cite{PeRu} a pioneer study on automatic continuity
of ternary derivations from a JB$^*$-triple $E$ into a Jordan Banach triple $E$-module. The concept of
Jordan Banach triple module is introduced in the just quoted paper, where it is also established
that every triple derivation from a real or complex JB$^*$-triple into its dual space or into itself is automatically
continuous. It seems natural, at this stage, to consider generalised triple derivations
in the context of JB$^*$-triples, studying the automatic continuity of these mappings.\smallskip

Jordan triple modules over Jordan triples were introduced as appropriate
extensions of bimodules over associative algebras and Jordan modules over Jordan algebras (cf. \cite{PeRu}). The concrete
definition reads as follows: Let $E$ be a complex (resp. real) Jordan triple,  a \emph{Jordan
triple $E$-module}  (also called \emph{triple $E$-module}) is a
vector space $X$ equipped with three mappings $$\{.,.,.\}_1 :
X\times E\times E \to X, \quad \{.,.,.\}_2 : E\times X\times E \to
X$$ $$ \hbox{ and } \{.,.,.\}_3: E\times E\times X \to X$$
satisfying  the following axioms:
\begin{enumerate}[{$(JTM1)$}]
\item $\{ x,a,b \}_1$ is linear in $a$ and $x$ and conjugate
linear in $b$ (resp., trilinear), $\{ a,b,x \}_3$ is linear in $b$
and $x$ and conjugate linear in $a$ (resp., trilinear) and
$\{a,x,b\}_2$ is conjugate linear in $a,b,x$ (resp., trilinear)
\item  $\{ x,b,a \}_1 = \{ a,b,x \}_3$, and $\{ a,x,b \}_2  = \{
b,x,a \}_2$  for every $a,b\in E$ and $x\in X$. \item Denoting by
$\J ...$ any of the products $\{ .,.,. \}_1$, $\{ .,.,. \}_2$ and
$\{ .,.,. \}_3$, the identity $\J {a}{b}{\J cde} = \J{\J abc}de $
$- \J c{\J bad}e +\J cd{\J abe},$ holds whenever one of the
elements  $a,b,c,d,e$ is in $X$ and the rest are in $E$.
\end{enumerate}

When $E$ is a Jordan Banach triple and $X$ is a triple $E$-module
which is also a Banach space, we shall say that $X$ is a \emph{Banach (Jordan) triple
$E$-module} when the products $\{ .,.,.\}_1$, $\{ .,.,. \}_2$ and $\{ .,.,. \}_3$ are (jointly)
continuous. From now on, the products $\{ .,.,.\}_1$, $\{ .,.,. \}_2$ and $\{ .,.,. \}_3$
will be simply denoted by $\J ...$.\smallskip

Every real or complex associative algebra $A$ (resp.,  Jordan
algebra $J$) is a real Jordan  triple with respect to $\J abc :=
\frac12 \left(abc +cba\right)$, $a,b,c\in A$ (resp., $\J abc =
(a\circ b) \circ c + (c\circ b) \circ a - (a\circ c) \circ b$) ,
$a,b,c\in J$). It is not hard to see that every $A$-bimodule $X$
is a real triple $A$-module with respect to the products $\J abx_3
:=   \frac12 \left(abx +xba\right )$ and  $\J axb_2 = \frac12
\left(axb +bxa\right)$, and that every Jordan module $X$ over a
Jordan algebra $J$ is a real triple $J$-module with respect to the
products $\J abx_3 :=  (a\circ b) \circ x + (x\circ b) \circ a -
(a\circ x) \circ b$ and  $\J axb_2 = (a\circ x) \circ b + (b\circ
x) \circ a - (a\circ b) \circ x.$\smallskip

The dual space, $E^*$, of a complex (resp., real) Jordan Banach triple $E$ is a
complex (resp., real) triple $E$-module with respect to the
products: \begin{equation}\label{eq module product dual 1}   \J
ab{\varphi} (x) = \J {\varphi}ba (x) := \varphi \J bax
\end{equation} and \begin{equation}\label{eq module product dual
2} \J a{\varphi}b (x) := \overline{ \varphi \J axb }, \end{equation}
$ \forall \varphi\in E^*, a,b,x\in E$ (cf. \cite{PeRu}).\smallskip

Given a triple $E$-module $X$ over a Jordan triple $E$, the
space $E\oplus X$ can be equipped with a structure of real Jordan
triple with respect to the product $\J {a_1+x_1}{a_2+x_2}{a_3+x_3}
= \J {a_1}{a_2}{a_3} +\J  {x_1}{a_2}{a_3}+\J  {a_1}{x_2}{a_3} + \J
{a_1}{a_2}{x_3}$. The Jordan triple $E\oplus X$ will be called the
\emph{triple split null extension} of $E$ and $X$.\smallskip

Let $X$ be a Jordan triple $E$-module over a Jordan triple $E$.
A \emph{triple derivation} from $E$ to $X$ is a conjugate linear map $\delta: E \to X$
satisfying $\delta \J abc = \J {\delta(a)}{b}{c} + \J {a}{\delta(b)}{c} + \J {a}{b}{\delta(c)}.$\smallskip

Let $E$ be a real (resp., complex) Jordan Banach triple and let $X$ be a Jordan Banach
triple $E$-module. A (conjugate) linear mapping $\delta:E\to X$ is said to be a \emph{generalised derivation}
when the mapping $\check{\delta} : E \times E \times E \to X,$
$$(a,b,c)\mapsto \check{\delta}(a,b,c):=\delta \{a,b,c\}-\{\delta(a),b,c\}-\{a,\delta(b),c\}-\{a,b,\delta(c)\}$$
is (jointly) continuous.\smallskip

Arguing as in \cite{PeRu}, we will associate to each generalised derivation from a JB$^*$-triple $E$
into a Jordan Banach triple $E$-module a generalised triple
homomorphism, in such a a way that the continuity of these two mappings is mutually determined.\smallskip

Let $\delta:E\to X$ be a generalised derivation. The symbol $E\oplus X$ will stand for the triple
split null extension of $E$ and $X$ equipped with the $\ell_{1}$-norm.
We define the mapping $$\Theta_\delta : E\to E\oplus X,$$ $$a\mapsto a +\delta(a).$$
It is clear that $\delta$ is continuous if and only if $\Theta_{\delta}$ is continuous. Furthermore,
the identity $$ \check{\delta}(a,b,c)=\delta \{a,b,c\}-\{\delta(a),b,c\}-\{a,\delta(b),c\}-\{a,b,\delta(c)\}$$
$$={\Theta}_\delta \J abc - \J {{\Theta}_\delta(a)}{{\Theta}_\delta(b)}{{\Theta}_\delta(c)} =\check{\Theta}_\delta(a,b,c), \ \forall a,b,c \in E,$$
shows that $\Theta_\delta$ is a generalised triple homomorphism. According to this notation, we set
$\Delta:=\Theta_{\delta}(E)=\{a+\delta(a):a\in  E\}.$ Let $(E\oplus X)_\Delta$ be the norm closed subtriple of $E\oplus X$
generated by $\Delta.$ Since $\Theta_{\delta}$ is a generalised triple homomorphism, by Lemma \ref{l separating spaces are
  ideals}, the separating space $\sigma_{_{E\oplus X}}(\Theta_{\delta})$ is a triple ideal of $(E\oplus X)_\Delta.$
It is not hard to see that $\sigma_{_{E\oplus X}}(\Theta_{\delta})$ coincides with $\{0\} \times \sigma_X(\delta).$\smallskip

A subspace $S$ of a triple $E$-module $X$ is said to be a
\emph{Jordan triple submodule} or a \emph{triple submodule} if and
only if $\J EES \subseteq S$ and $ \J ESE \subseteq S$. Every triple ideal $J$ of
$E$ is a Jordan triple $E$-submodule of $E$.\smallskip

Let $a+x,b+y$ be elements in $(E\oplus X)_\Delta$ and $z\in \{0\} \times  \sigma_X(\delta)=\sigma_{_{E\oplus X}}(\Theta_{\delta}).$
By the definition of the triple product in $E\oplus X$ and the just quoted fact that $\sigma_{_{E\oplus X}}(\Theta_{\delta})$
is a triple ideal of $(E\oplus X)_\Delta$ we have \begin{equation}\label{eq 1 gender}
\{a,b,z\}=\{a+x,b+y,z\} \in \sigma_X(\delta)
  \end{equation}  and \begin{equation}\label{eq 2 gender}
\{a,z,b\}=\{a+x,z,b+y\} \in \sigma_X(\delta).
  \end{equation}

Since $(E\oplus X)_\Delta$ contains $\Delta$, it follows from $(\ref{eq 1 gender})$ and $(\ref{eq 2 gender})$ that
$\{E,E,\sigma_X(\delta)\} \subseteq \sigma_X(\delta)$ and $\{E,\sigma_X(\delta),E\} \subseteq \sigma_X(\delta).$
Since $\sigma_X(\delta)$ is always a linear subspace, it also is a triple $E$-submodule of $X.$\smallskip

For each subset $A$ of a triple $E$-module $X$, we define its
\emph{quadratic annihilator}, Ann$_{E} (A)$, as the set
$\{ a\in E : Q (a) (A) = \J aAa = 0\}$.\smallskip

We shall also make use of the following equality:
\begin{equation}\label{eq 3 gender}
Ann_{E\oplus X}(\sigma_{E\oplus X}(\Theta_\delta))=Ann_E(\sigma_X(\delta))\oplus X.
\end{equation}

\begin{remark}\label{r annih E,E*}{\rm The quadratic annihilator of
a submodule $S$ of a triple module $X$ need not be, in general, a
linear subspace (cf. \cite{PeRu}). However, it is known that when $E$ is a JB$^*$-triple
and $X=E$ or $X=E^*$ then, for each submodule $S$ of $X$, $Ann_E(S)$
is a linear subspace, and hence a norm-closed triple ideal of $E$
(see Lemma 1 and Proposition 2 in \cite{PeRu}). Further, Proposition 2 (or Remark 3)
in \cite{PeRu} shows that, in this case,
$\{Ann_E(S),Ann_E(S),S\}=0$ in the triple split null extension $E\oplus X.$}
\end{remark}

From now on, we assume that $E$ is a JB$^*$-triple and $X$ denotes $E$ or $E^*.$
In this case, Remark \ref{r annih E,E*} and the fact that $\sigma_X(\delta)$ is a triple
$E$-submodule of $X$ prove that $Ann_E(\sigma_X(\delta))$ is a norm-closed triple ideal of $E.$\smallskip

The strategy to obtain results on automatic continuity for generalised triple derivations will consist
in applying Theorem \ref{t gen.caract.cont} to the generalised triple homomorphism $\Theta_{\delta}.$
In order to do this, we shall first check that
$$J:=\Theta_{\delta}^{-1}(Ann_{E\oplus X}(\sigma_{E\oplus
X}(\Theta_{\delta})))$$ is a norm-closed triple ideal of $E.$ It is
not hard to see that $Ann_{E\oplus X}(\sigma_{E\oplus
X}(\Theta_{\delta})))=Ann_E(\sigma_X(\delta))\oplus X$ and
$$\Theta_{\delta}^{-1}(Ann_E(\sigma_X(\delta))\oplus
X)=Ann_E(\sigma_X(\delta)).$$

This proves that $J$ is a norm-closed triple ideal of $E$ (see
remark \ref{r annih E,E*}). On the other hand,
$$\{ Ann_{E\oplus X}( \sigma_{E\oplus X}(\Theta_{\delta}), Ann_{E\oplus X}( \sigma_{E\oplus
 X}(\Theta_{\delta})), \sigma_{E\oplus X}(\Theta_{\delta})\}$$ $$=\{Ann_E(\sigma_X(\delta)),Ann_E(\sigma_X(\delta)),\sigma_X(\delta)\}=0$$
(compare the final statement in Remark \ref{r annih E,E*}).  Theorem \ref{t gen.caract.cont} proves the continuity of $\Theta_{\delta}$
and hence the continuity of $\delta$.

\begin{theorem}\label{t genderJB*-triple} Let $E$ be a real or complex JB$^*$-triple and $\delta:E\to X$ a
generalised triple derivation, where $X=E$ or $E^*$. Then $\delta$ is continuous.$\hfill \Box$
\end{theorem}

The statement concerning real JB$^*$-triples can be derived from the complex case applying
Remark 14 in \cite{PeRu}.\smallskip

Since every triple derivation is a generalised triple derivation we get the following:

\begin{corollary}\label{t DerJB*-triple}\cite[Corollary 15]{PeRu} Let $E$ be a real or complex
JB$^*$-triple and let $\delta:E\to X$ be a triple derivation, where $X=E$
or $E^*$. Then $\delta$ is continuous. $\hfill \Box$
\end{corollary}

\subsection{Generalised triple derivations whose domain is a C$^*$-algebra}

We have already mentioned that every C$^*$-algebra belongs to the class of JB$^*$-triples.
We shall conclude this paper by applying some of the previous results to
C$^*$-algebras. The results obtained this way are interesting by themselves.

\begin{lemma}\label{l techn.CstarToBanach}
Let $T:A_{sa}\to X$ be a linear mapping from the self-adjoint part, $A_{sa}$,
of an abelian C$^*$-algebra, $A$, to a Banach space.
Suppose that $J_T:=\{ a\in A_{sa} : T Q(a) \hbox{ is continuous }\}$
is a norm-closed closed subset of $A_{sa}$ with $\{a,A_{sa},a\}\in J_T,$
for every $a\in J_{T}$. Then $J_T$ is a triple ideal of $A_{sa}.$
\end{lemma}

\begin{proof} It is easy to see that every norm-closed inner ideal of the selfadjoint
part of an abelian C$^*$-algebra $A$ is a norm-closed triple ideal of $A_{sa}$.
Therefore, we only have to prove that $J_T$ is a linear subspace, to this end, it is enough to
show that $a+b\in J_{T}$ whenever $a,b\in J_T$.\smallskip

Let $a,b$ be two elements in $A_{sa}.$ First we observe that, since
$A_{sa}$ is abelian, $L(a+b)=Q(a+b)$. Obviously, the linear mapping $L_b:A_{sa}\to A_{sa}, c\mapsto
cb = bc$ is continuous. Since $A_{sa}$ is abelian we have
$L(a^2,b)=Q(a) L_b= L_b Q(a).$ Therefore $TL(a^2,b)=TQ(a)L_b$ is
continuous, for every $a\in J_{T}$, $b\in A_{sa}$.\smallskip

Let us pick $a\in J_{T}$.  We write $a$ in the
form $a=a_1-a_2$ where $a_1,a_2$ are orthogonal positive elements
in $A_{sa}.$ Since $Q(a)A_{sa} \in J_T$, $a_1^3$ lies in $J_T$,
and hence $a_1^6 \ A_{sa}=Q(a_1^3) A_{sa} \subseteq J_T.$
This implies that $J_T$ contains the norm-closed ideal of $A_{sa}$ generated by $a_1^6$,
which guarantees that $J_T$ contains $a_1$ and $a_1^{\frac{1}{2}}.$ Similarly, we show $J_T$
contains $a_2$ and $a_2^{\frac{1}{2}}.$ Now
$$T L(a,b)=TL(a_1,b)-TL(a_2,b) = T L((a_1^{\frac{1}{2}})^2,b) - TL((a_2^{\frac{1}{2}})^2,b),$$
and thus $T L(a,b)$ is continuous for every $b\in A_{sa}$. Finally, the equality
$$TQ(a+b) = T L(a+b)=T L(a,a)+T L(b,b)+2 T L(a,b)$$
shows that $T Q(a+b)$ is continuous for every $a,b\in J_{T}$.
\end{proof}

\begin{proposition}\label{p abelianTobimodule}
Let $\delta:A\to X$ be a generalised derivation from an abelian
C$^*$-algebra to a Jordan Banach triple $A$-module. Then $\delta$ is
continuous.
\end{proposition}

\begin{proof} We shall only prove that $\delta_{|A_{sa}}$ is continuous.
Let $\Theta_{\delta_{0}}: A_{sa} \to A_{sa}\oplus X$ be the
generalised triple homomorphism associated to $\delta_{0}:=\delta|_{A_{sa}}.$
We have already shown that $J=\Theta_{\delta_{0}}^{-1} (Ann_{A_{sa}\oplus
X}(\sigma_{A_{sa}\oplus X}(\Theta_{\delta_{0}})))$ coincides
with $Ann_{A_{sa}}(\sigma_{X}(\delta_{0}))$ (see the comments
prior to Theorem \ref{t genderJB*-triple}). Therefore, $J$ is the
quadratic annihilator of a closed submodule of $X$, and hence $J$ is norm closed
and satisfies $\{a,A_{sa},a\}\in J,$ for every $a\in J$ (cf. \cite[\S 2.3]{PeRu}).\smallskip

It is easy to see that $J$ coincides with $\{a\in A_{sa} : \Theta_{\delta_{0}} Q(a) \hbox{ is continuous } \}.$
Now, Lemma \ref{l techn.CstarToBanach} proves that $J=\Theta_{\delta_{0}}^{-1} (Ann_{A_{sa}\oplus
X}(\sigma_{A_{sa}\oplus X}(\Theta_{\delta_{0}})))$ is a norm-closed triple ideal of $A_{sa}$, and since $A$
is abelian $$\J {Ann_{A_{sa}\oplus X}(\sigma_{A_{sa}\oplus X}(\Theta_{\delta_{0}}))}{Ann_{A_{sa}\oplus
X}(\sigma_{A_{sa}\oplus X}(\Theta_{\delta_{0}}))}{\sigma_{A_{sa}\oplus X}(\Theta_{\delta_{0}})}$$ $$=\J {Ann_{A_{sa}}(\sigma_{A_{sa}}(\delta_{0}))}{Ann_{A_{sa}}(\sigma_{A_{sa}}(\delta_{0}))}{\sigma_{A_{sa}}(\delta_{0})} =0.$$

Having in mind that $A_{sa}$ is a reduced real JB$^*$-triple and the validity of Theorem \ref{t gen.caract.cont}
for reduced real JB$^*$-triples, we prove that $\delta|_{A_{sa}}$ is continuous.\smallskip
\end{proof}

A celebrated result of J. Cuntz (see \cite{Cuntz}) establishes that a linear mapping from a C$^*$-algebra $A$
to a Banach space is continuous if and only if its restriction to each subalgebra of $A$ generated
by a single hermitian element is continuous. We finish this note with a consequence of Cuntz theorem
and Proposition \ref{p abelianTobimodule}.

\begin{theorem}\label{t CtsarToBimodule}
Every generalised triple derivation from a real or complex C$^*$-algebra $A$ to a
Jordan Banach triple $A$-module is continuous. $\hfill \Box$
\end{theorem}

\bigskip\bigskip

\end{document}